\documentclass[11pt]{amsart}
\usepackage{amssymb}
\usepackage{amsthm}

\usepackage{graphicx}
\usepackage{epstopdf}
\usepackage{verbatim}

\def\Z{{\mathbb Z}}

\def\a{{\alpha}}
\def\b{{\beta}}
\def\g{{\gamma}}
\def\d{{\delta}}

\def\TSG{{\mathrm{TSG_+}}}
\def\Aut{{\mathrm{Aut}}}
\def\Diff{{\mathrm{Diff_+}}}

\def\cl{{\mathrm{cl}}}

\def\kernel{{\mathrm{\ker}}}
\newcommand{\x}{\times}
\newcommand{\bd}{\partial}

\newcommand{\rank}{\mathrm{rank}}

\newcommand{\sbgp}{\leq}

\newtheorem*{Subgroup Corollary}{Subgroup Corollary}
\newtheorem*{Subgroup thm}{Subgroup Theorem}

\newtheorem*{add}{Knot Addition Lemma}
\newtheorem*{auto}{Automorphism Theorem}
\newtheorem*{complete}{Complete Graph Theorem}
\newtheorem*{finite}{Finiteness Theorem}

\newtheorem*{Z3}{{\bf $\mathbb{Z}_3\times\mathbb{Z}_3$} Lemma}

\newtheorem{prop}{Proposition}

\newtheorem{definition}{Definition}
\newtheorem{thm}{Theorem}

\newtheorem*{JSJ}{Characteristic Submanifold Theorem [JSJ]}

\newtheorem{lemma}{Lemma}

\title[Spatial graphs with local knots]{Spatial graphs with local knots}
\author{Erica Flapan, Blake Mellor, and Ramin Naimi}

\thanks{The authors were supported in part 
by NSF  
 Grants DMS-0905087, DMS-0905687, DMS-0905300.}



\address{Department of Mathematics, Pomona College, Claremont, CA 91711, USA}

\address{Department of Mathematics, Loyola Marymount University, Los Angeles, CA 90045, USA}

\address{Department of Mathematics, Occidental College, Los Angeles, CA 90041, USA}

\begin{document}
\date \today
\maketitle
\begin{abstract}
It is shown that for any locally knotted edge of a 3-connected graph in $S^3$, there is a ball that contains all of the local knots of that edge which is unique up to an isotopy setwise fixing the graph.  This result is applied to the study of topological symmetry groups of graphs embedded in $S^3$.
\end{abstract}
\bigskip

Schubert's  1949 result \cite{Sch} that every non-trivial knot can be uniquely factored into prime knots is a fundamental result in knot theory.   Hashizume \cite{Ha}, extended Schubert's result to links in 1958.  Then in 1987, Suzuki \cite{Su} generalized Schubert's result to spatial graphs by proving that every connected graph embedded in $S^3$ can be split along spheres meeting the graph in 1 or 2 points to obtain a unique collection of prime embedded graphs together with some trivial graphs.   

\begin{figure} [h]
\includegraphics{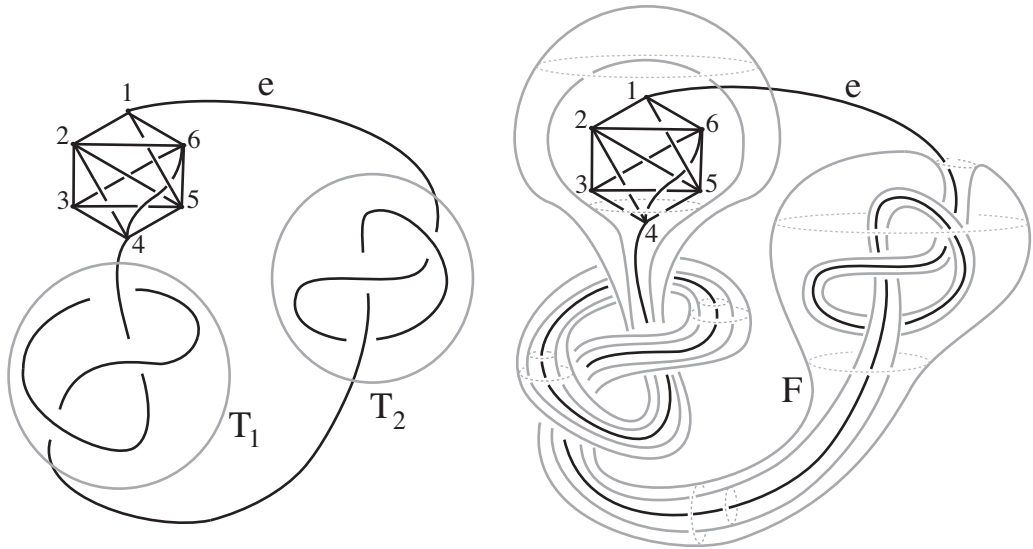}
\caption{The sphere $F$ on the right is not isotopic (setwise fixing the graph) to either of the spheres $T_1$ or $T_2$ on the left.}
\label{nonstandard}
\end{figure}

Although the set of prime factors of a knot or embedded graph is unique up to equivalence, the set of splitting spheres is generally not unique up to an isotopy setwise fixing the knot or graph. For example, consider the embedding of the complete graph $K_6$ which is illustrated on both the left and right sides of Figure \ref{nonstandard}.  The edge $e=\overline{14}$ contains two trefoil knots.  The spheres $T_1$ and $T_2$ (illustrated on the left) are splitting spheres for these two knots.  However, one of the balls bounded by $F$ (illustrated on the right) meets $e$ in an arc whose union with an arc in $F$ is a single trefoil knot.  Thus $F$ is also a splitting sphere for one of the two local knots in $e$.  However, $F$ is not isotopic (fixing the embedded graph setwise) to either of the spheres $T_1$ or $T_2$.

By contrast, in this paper we show that for any locally knotted edge of an embedded 3-connected graph, there is a ball meeting the graph in an arc containing all of the local knots of that edge which is unique up to an isotopy fixing the graph. We call such a ball an {\it unknotting ball} for that edge. Our main theorem is the following.
\medskip

\begin{thm}\label{P:uniqueness}Let $\Gamma$ be a 3-connected graph embedded in $S^3$.  Then any locally knotted edge $e$ has an unknotting ball which is unique up to an isotopy of $(S^3,\Gamma)$ fixing the vertices of $\Gamma$. Furthermore, if $B_{1},\dots,B_{n}$ are pairwise disjoint balls for local knots of an edge $e$, then $e$ has an unknotting ball which contains $B_{1} \cup \cdots \cup B_{n}$.

\end{thm}
\medskip

We shall apply this theorem to the study of topological symmetry groups of graphs embedded in $S^3$. The concept of the {\it topological symmetry group} 
was first introduced by Jon Simon \cite{Si} as a way of describing the
symmetries of non-rigid molecules.  Let $\gamma$ be an abstract graph.  An {\it automorphism} of $\gamma$ is a permutation of the vertices of $\gamma$ which preserves adjacency.  We use $\Aut(\gamma)$ to denote the group of automorphisms of $\gamma$.  Given an embedding $\Gamma$ of  an abstract graph $\gamma$ in $S^3$, the {\it topological symmetry group}, ${\rm TSG}(\Gamma)$, is defined to be the subgroup of $\Aut(\gamma)$ induced on the vertices of $\gamma$ by diffeomorphisms of the pair $(S^3, \Gamma)$.  If we only allow orientation preserving diffeomorphisms, we obtain the orientation preserving topological symmetry group $\TSG(\Gamma)$.   In this paper we are only interested in the orientation preserving topological symmetry group.  Thus for simplicity, we abuse notation and refer to $\TSG(\Gamma)$ as the {\it topological symmetry group} rather than the {\it orientation preserving topological symmetry group}.  

Flapan, Naimi, Pommersheim, and Tamvakis \cite{FNPT} proved that not every finite group can occur as  $\TSG(\Gamma)$ for some embedded graph $\Gamma$ in $S^3$.  For example, the alternating groups $A_n$ for $n>5$ cannot occur as $\TSG(\Gamma)$ for any embedded graph.  For most abstract graphs $\gamma$, it is not known what groups can occur as $\TSG(\Gamma)$ for some embedding of $\gamma$ in $S^3$.  
In this paper, we consider the simpler question of whether all of the subgroups of a given $\TSG(\Gamma)$ can occur as $\TSG(\Gamma')$ for some re-embedding $\Gamma'$ of $\Gamma$ in $S^3$.  

\begin{figure} [h]
\includegraphics{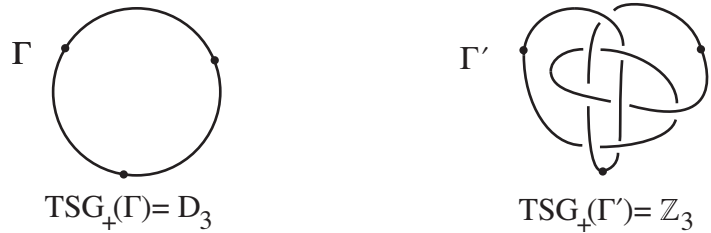}
\caption{Embeddings of a triangle with different topological symmetry groups.}
\label{D3}
\end{figure}

For example, let $\gamma$ denote a circle with three vertices.   If the embedded graph $\Gamma$ is an unknotted circle, then $\TSG(\Gamma)=D_3$, the dihedral group with 6 elements.  If we re-embed $\Gamma$ as a non-invertible knot $\Gamma'$, then $\TSG(\Gamma')$ is the subgroup $\mathbb{Z}_3$ (in Figure \ref{D3}, $\Gamma'$ contains the non-invertible knot $8_{17}$).  On the other hand, for any embedding $\Gamma''$ of $\gamma$ in $S^3$ there will be an orientation preserving diffeomorphism of $(S^3, \Gamma'')$ which cyclically permutes the three vertices (obtained by slithering $\Gamma''$ along itself) inducing an order 3 automorphism of $\Gamma''$.  Hence there is no embedding $\Gamma'''$ of $\gamma$ such that $\TSG(\Gamma''')$ is either $\mathbb{Z}_2$ or the trivial group.  Thus not all of the subgroups of $\TSG(\Gamma)$ can occur as $\TSG(\Gamma''')$ for some re-embedding $\Gamma'''$ of $\Gamma$ in $S^3$.  

In contrast with the above example, we use local knotting together with Theorem~1 to prove that in many cases  every subgroup of a topological symmetry can occur as the topological symmetry group of another embedding of the graph.  In particular, we prove the following.

\medskip

\begin{Subgroup thm}  Let $\Gamma$ be an embedding of a 3-connected graph in $S^3$ and let $H$ be a (possibly trivial) subgroup of $\TSG(\Gamma)$.  Let $e_1$, \dots $e_n$ be a set of edges in $\Gamma$ whose orbits under $H$ are distinct.  Suppose that any $\varphi\in \TSG(\Gamma)$ which  pointwise fixes $e_1$ and satisfies $\varphi(\langle e_i\rangle_H)=\langle e_i\rangle_H$ for each $i$, also pointwise fixes a subgraph of $\Gamma$ that cannot be embedded in $S^1$.  Then there is an embedding $\Gamma'$ of $\Gamma$ with $H=\TSG(\Gamma')$.
\end{Subgroup thm}
\medskip

  Observe that the result below follows immediately from the Subgroup Theorem.

\medskip
\begin{Subgroup Corollary} \label{L:fundamentaledge}
Let $\Gamma$ be a 3-connected graph embedded in $S^3$ which has an edge $e$ that is not pointwise fixed by any non-trivial element of $\TSG(\Gamma)$.  Then for every (possibly trivial) subgroup $H$ of $\TSG(\Gamma)$ there is an embedding $\Gamma'$ of $\Gamma$ with $H = \TSG(\Gamma')$.
\end{Subgroup Corollary}

\medskip

We shall apply this corollary to the study of the topological symmetry groups of complete graphs.  The complete graphs $K_n$ are an interesting family of graphs to focus on because the automorphism group of $K_n$ is the symmetric group $S_n$, which is the largest automorphism group of any graph
with $n$ vertices. On the other hand, Flapan, Naimi, and Tamvakis \cite{FNT} have classified the groups that can occur as the topological symmetry group of some embedding of a complete graph in $S^3$. More specifically, they proved the following.

\medskip

\begin{complete} \label{T:FNT} \cite{FNT}
A finite group $H$ is isomorphic to $\TSG(\Gamma)$ for some embedding $\Gamma$ of a complete graph in $S^3$ if and only if $H$ is a finite cyclic group, a dihedral group, a subgroup of $D_m \times D_m$ for some odd $m$, or $A_4$, $S_4$, or $A_5$.
\end{complete}
\medskip

Although this result restricts the types of groups which can occur, for a given complete graph $K_n$, it is still not known precisely which of the above groups occur.  In this paper, we use the Subgroup Corollary to prove the following.  

\medskip
\setcounter{thm}{2}
\begin{thm} \label{T:realize}
Let $n > 6$ and let $\Gamma$ be an embedding of $K_n$ in $S^3$ such that $\TSG(\Gamma)$ is a finite cyclic group, a dihedral group, or a subgroup of $D_m \times D_m$ for some odd $m$.  Then for every (possibly trivial) subgroup $H$ of $\TSG(\Gamma)$, there is an embedding $\Gamma'$ of $K_n$ such that $H = \TSG(\Gamma')$.
\end{thm}

\medskip

Note that in \cite{FMN}, we classify all values of $n$ such that there is an embedding $\Gamma$ of $K_n$ with $\TSG(\Gamma)$ equal to $A_4$, $S_4$, or $A_5$.

 \bigskip

 \section{Unknotting balls}

We begin with some terminology.   By a {\it graph} we mean a finite set of vertices and edges such that there is at most one edge between a pair of vertices and every edge has two distinct vertices.  A graph is said to be  {\it 3-connected} if it cannot be disconnected or reduced to a single vertex by removing fewer than 3 vertices together with the edges containing them.  For example, the complete graph $K_n$ is 3-connected if and only if $n>3$.  

Let $\Gamma$ be a 3-connected graph embedded in $S^3$ and let $e$ be an edge of $\Gamma$.  If there is a ball $B$ in $S^3$ such that $B\cap \Gamma$ is an arc in the interior of $e$ whose union with an arc in $\partial B$ has non-trivial knot type $K$, then we say that $B$ is a {\it ball for the local knot $K$ of $e$} and we say that the pair $(B,B\cap \Gamma)$ has {\it knot type $K$}.  It is shown in \cite{FNPT} that for any edge of an embedded 3-connected graph the local knots on that edge are well defined.  If a graph is not 3-connected this is not necessarily the case.  For example, see Figure 3. 

\begin{figure} [h]
\includegraphics{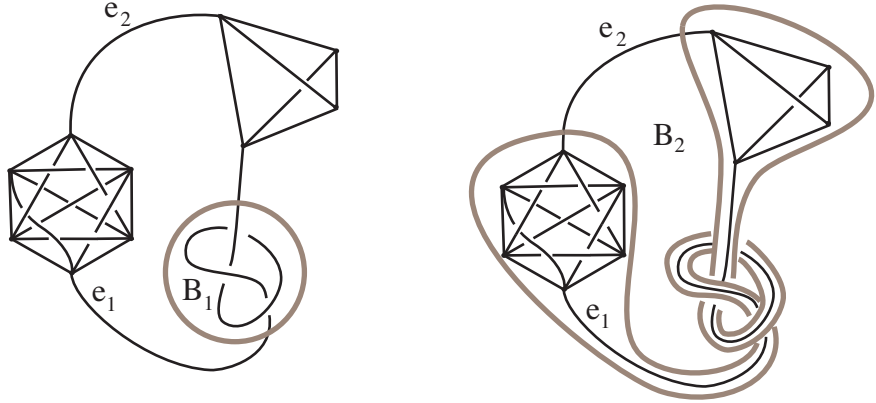}
\caption{On the left we see a local knot on $e_1$ with ball $B_1$, and on the right we see the ``same knot'' is a local knot on $e_2$ with ball $B_2$.}
\label{decomp}
\end{figure}

 Suppose that some ball $B$ meets $\Gamma$ in an arc $A=B\cap \Gamma$.  If $(B,A)$ has non-trivial knot type then we say the arc $A$ is {\it knotted}, otherwise we say the arc $A$ is {\it unknotted}.  Let $e$ be an edge of $\Gamma$ which contains a local knot with ball $B$.   Let $e'$ be an arc obtained from $e$ by replacing the knotted arc $e\cap B$ by an unknotted arc of $B$ with the same endpoints and let $\Gamma'$ be obtained from $\Gamma$ by replacing $e$ by $e'$.  If $e'$ has no local knots in $\Gamma'$, then we say that $B$ is an {\it unknotting ball} for $e$.  We show below that every locally knotted edge has an unknotting ball and this ball is unique up to an isotopy setwise fixing $\Gamma$.  This is in contrast with the observation that balls for local knots are generally not unique up to isotopy as illustrated in Figure 1.

 We  construct a neighborhood of $\Gamma$ out of balls and tubes  around vertices and edges as follows.  Let $V$ and $E$ denote the vertices and edges respectively of $\Gamma$.  For each vertex $v\in V$, we let $N(v)$ 
denote a closed ball around $v$ whose intersection with $\Gamma$ consists of the single vertex $v$ together with an arc of every edge containing $v$ such that if $w\not =v$ then $N(w)\cap N(v)=\emptyset$.  Let $N(V)$ 
denote the union of all of these balls.  For each 
embedded edge $e \in E$, let 
$N(e )$ denote a solid cylinder $D^2\times I$
whose core is $e -N(V)$, such that 
$N(e)\cap \Gamma\subseteq e$, $N(e)\cap N(e')=\emptyset$ for $e'\not =e$, and $N(e)$ meets $N(V)$ in 
a pair of disks $D\times \{0,1\}$.  Let $N(E)$ 
denote the union of all these solid cylinders, and let $N(\Gamma )=N(V)\cup N(E)$.  Suppose that $e$ is an edge in an embedded graph $\Gamma$ with vertices $v_1$ and $v_2$, and $B$ is a ball such that $B\cap \Gamma$ is an arc of $e$.  Observe that there is an isotopy of $(S^3,\Gamma)$ fixing every vertex of $\Gamma$ taking $B$ to a ball whose boundary meets $N(\Gamma)$ in two disks contained in $\partial N(v_1)$ and $\partial N(v_2)$.  Thus we shall assume, when needed, that a ball for a local knot of $e$ has this form.

\medskip

We will use the Characteristic Submanifold Theorem for pared manifolds together with some definitions stated below.  

\medskip\nobreak 

\begin{definition}    A {\bf pared 
$3$-manifold} $(M,P)$ is an orientable 
$3$-manifold $M$ together with a family $P$ of 
disjoint incompressible annuli and tori in 
$\partial M$.\end{definition}      

\medskip\nobreak

\begin{definition} A pared manifold 
$(M,P)$ is said to be {\bf simple} if it 
satisfies the following three conditions: 
     
\par \noindent 1)  $M$ is irreducible and 
$\partial M-P$ is incompressible

\par \noindent 2)  Every incompressible torus in 
$M$ is parallel to a torus component of $P$

\par \noindent 3)  Any annulus $A$ in $M$ with 
$\partial A$ contained in $\partial M-P$ is 
either compressible or parallel to an annulus 
$A'$ in $\partial M$ with $\partial A'=\partial 
A$ and such that $A'\cap P$ consists of zero or 
one annular component of $P$.\end{definition}    

\medskip\nobreak 

\begin{definition}     A pared manifold 
$(M,P)$ is said to be {\bf Seifert fibered} if 
there is a Seifert fibration of $M$ for which $P$ 
is a union of fibers.  A pared manifold $(M,P)$ 
is said to be {\bf I-fibered} if there is an 
$I$-bundle map of $M$ over a surface $S$ such 
that $P$ is the preimage of $\partial 
S$. \end{definition}      

\medskip\nobreak

We use the following version of the Characteristic Submanifold Theorem for pared manifolds proved independently by Jaco and Shalen \cite{JS} and Johannson \cite{Jo}.  Henceforth we shall refer to this result simply as JSJ.

\medskip\nobreak 

 \begin{JSJ} \cite{JS,Jo} Let $(M,P)$ be a pared manifold with 
$M$ irreducible and $\partial M-P$ 
incompressible.  Then, up to an isotopy of 
$(M,P)$, $M$ contains a unique finite family $\tau$ of 
disjoint incompressible tori and annuli with boundaries in
$\partial M-P$ 
with the following two defining properties:
\medskip

\par \noindent 1)  If $Q$ is the closure of a 
component of $M-\tau $, then the pared manifold 
$(Q,Q\cap (P\cup \tau ))$ is either simple, 
Seifert fibered, or $I$-fibered. 
\medskip

\par \noindent 2)  No such family has
fewer elements than $\tau$.
\medskip

\noindent Furthermore, if $F$ is an annulus in $M$ with boundaries in $\partial M-P$ which is incompressible and boundary incompressible, then there is an isotopy of the pair $(M,P\cup \tau)$ taking $F$ to an annulus which is contained in either $\tau $ or in a component $(Q,Q\cap (P\cup \tau ))$ that is
Seifert fibered or $I$-fibered. 

\end{JSJ}    
\medskip

We now state and prove Theorem~1.

\medskip

\setcounter{thm}{0}

\begin{thm}\label{P:uniqueness}
Let $\Gamma$ be a 3-connected graph embedded in $S^3$.  Then any locally knotted edge $e$ has an unknotting ball which is unique up to an isotopy of $(S^3,\Gamma)$ fixing the vertices of $\Gamma$. Furthermore, if $B_{1},\dots,B_{n}$ are pairwise disjoint balls for local knots of an edge $e$, then $e$ has an unknotting ball which contains $B_{1} \cup \cdots \cup B_{n}$.

\end{thm}

\medskip\begin{proof}  Let $M=\cl(S^3-N(\Gamma))$ and $P=\cl(\partial N(E)-\partial N(V))$, and let $\tau$ be a JSJ family for $(M,P)$.  Thus each annulus in $P'=P\cup \tau$ has its boundary in $\partial N(V)$.   We begin by making two observations which will be used later.

First suppose, for the sake of contradiction, that there is a component of $M-\tau$ whose closure $(Q, Q\cap P')$ is $I$-fibered.  Since $M$ is contained in $S^3$, $Q$ is a product of a surface cross an interval where the ends of the product are in one or two components of $\partial N(V)$. Since $\Gamma$ is 3-connected, removing the vertices corresponding to these components and the edges containing them does not separate $\Gamma$.   Also, there is at most one edge between any pair of vertices and no edge from a vertex to itself.   It follows that there are precisely two annuli in $Q\cap P'$, and hence these annuli must be parallel in the product.  Since this contradicts the minimality of $\tau$, no components of $M-\tau$ are $I$-fibered.

  Next consider an annulus in $\tau$ whose boundary components are in a single $\partial N(v)$, and suppose, for the sake of contradiction, that its boundaries are non-isotopic curves in $\partial N(v)\cap M$.  Thus the surface between these two curves in $\partial N(v)$ intersects $\Gamma$.  Also, since the annulus is incompressible in $ M$, the disjoint disks bounded by the curves in $\partial N(v)$ each intersect $\Gamma$.  It follows that the sphere obtained by capping off the annulus with these two disks separates $\Gamma$ and each component contains at least one vertex other than $v$.  Since this contradicts the 3-connectedness of $\Gamma$, the boundary components of such an annulus must be isotopic in $\partial N(v)\cap M$.

We now consider the locally knotted edge $e$.  Let $v_1$ and $v_2$ be the vertices of $e$ and choose a ball for a local knot of $e$ whose boundary meets $N(\Gamma)$ in the disks $\partial N(v_1)\cap N(e)$ and $\partial N(v_2)\cap N(e)$.   By removing these disks we obtain an annulus which is incompressible in $M$.  It now follows from JSJ that there is an isotopy of the pair $(M,P')$ taking this annulus to an annulus $A$ which has boundaries in $\partial N(v_1)$ and $\partial N(v_2)$ and is contained either in $\tau$ or in a Seifert fibered component of $M-\tau$.  

Consider the case where the annulus $A$ is contained in a Seifert fibered component $(Q,Q\cap P')$ of $M-\tau$.   Let $D$ be the disk in $\partial N(v_1)-e$ bounded by $A\cap \partial N(v_1)$.   It follows from the fibered structure of $Q$ that every component of $Q\cap \partial N(v_1)$ is an annulus.  In particular, the component of $Q\cap \partial N(v_1)$ containing $\partial D$ has precisely one boundary which is in $D$ and is isotopic to $\partial D$ in $\partial N(v_1)\cap M$, and all of the other components of $Q\cap \partial N(v_1)$ have either 0 or 2 boundary components which are in $D$ and are isotopic to $\partial D$ in $\partial N(v_1)\cap M$.  Since every boundary component of $Q\cap \tau$  is a boundary of some component of $Q\cap \partial N(v_1)$, it follows that there are an odd number of boundary components of $Q\cap\tau$ which are in $D$ and are isotopic to $\partial D$ in $\partial N(v_1)\cap M$.  

Suppose, for the sake of contradiction, that every annulus of $Q\cap \tau$ with at least one boundary component which is in $D$ and is isotopic to $\partial D$ in $\partial N(v_1)\cap M$ has its other boundary component in $\partial N(v_1)$.   Such an annulus has isotopic boundaries in $\partial N(v_1)\cap M$, and since $A$ separates $Q$ and one boundary of it is in $D$ the other boundary must be in $D$ as well.   Thus there must be an even number of boundary components of $Q\cap\tau$ which are in $D$ and are isotopic to $\partial D$ in $\partial N(v_1)\cap M$.   Since we saw above that this number is actually odd, there must be some such annulus $A'$ whose boundary component is in $\partial N(v)$ with $v\not= v_1$.  

Now suppose, for the sake of contradiction, that $v\not =v_2$.  Consider the sphere obtained from $A\cup A'$ by adding the annulus between $A$ and $A'$ in $\partial N(v_1)\cap M$ together with disks in $N(v_2)$ and $N(v)$.  Since $A$ and $A'$ are incompressible in $M$, this sphere separates $\Gamma$ into components each containing at least one vertex of $\Gamma$ other than $v_2$ and $v$.  Since this contradicts the 3-connectedness of $\Gamma$, we must have $v=v_2$.  Furthermore, if we cap off $A'$ in $\partial N(v_1)$ and $\partial N(v_2)$, we obtain a sphere bounding a ball which intersects $\Gamma$ in an arc of $e$ and contains the annulus $A$.

Since $\Gamma$ is 3-connected, we can cap off any annulus in $\tau$ with boundaries in $\partial N(v_1)$ and $\partial N(v_2)$ to obtain a sphere bounding a ball whose intersection with $\Gamma$ is an arc of $e$.  Choose $F$ to be maximal among all such annuli.  That is, by capping off $F$ in $\partial N(v_1)$ and $\partial N(v_2)$ we obtain a sphere bounding a ball $B$ which intersects $\Gamma$ in an arc of $e$ and contains all of the other such annuli.  Then it follows from our argument above that any incompressible annulus in a Seifert fibered component of $M-\tau$ with boundaries in $\partial N(v_1)$ and $\partial N(v_2)$ will be contained in $B$.   Hence any incompressible annulus in $M$ with boundaries in $\partial N(V_1)$ and $\partial N(v_2)$ is isotopic in $(M,P')$ to an annulus which is properly embedded in $B\cap M$.  Furthermore, if such an annulus is contained in $\cl(M-B)$, then it is isotopic in $(\cl(M-B),P')$ to $F$.

It follows from our choice of $F$ that $B$ is an unknotting ball for $e$.  Let $B'$ be another unknotting ball for $e$,  and let $e'$ be obtained from $e$ by replacing the knotted arc $B'\cap e$ by an unknotted arc of $B'$. After an isotopy of $(S^3,\Gamma)$ fixing the vertices, we can assume that $\partial B'$ meets $N(\Gamma)$ in disks in $\partial N(v_1)$ and $\partial N(v_2)$.  By removing these disks we obtain an annulus $F'$, which we prove as follows is isotopic to $F$ in $(M,P')$. Since $B'$ is a ball for a local knot of $e$, $F'$ is incompressible in $M$.  Hence, by our choice of $F$, we can assume that $F'$ is properly embedded in $B\cap M$.  Let $T$ be the torus obtained from the annuli $F$ and $F'$ by adding annuli in $\partial N(v_1)$ and $\partial N(v_2)$.  Suppose that $F'$ is not isotopic to $F$.  If $T$ bounds a solid torus in $M$, then the meridian of the solid torus does not have intersection number $\pm1$ with a component of $\partial F$.  However, in this case, by adding a thickened disk to the solid torus along a component of $\partial F$ we would get a punctured lens space.  As this is impossible in $S^3$, we can assume that $T$ bounds a knot complement in $M$.  However, since $F'\subseteq B$, this would imply that $B$ is a ball for a local knot in the edge $e'$ which is contrary to our assumption that $B'$ is an unknotting ball for $e$.  Therefore $F'$ is isotopic to $F$ in $(M,P')$, and hence $B'$ is isotopic to $B$ by an isotopy of $(S^3, \Gamma)$ pointwise fixing the vertices of $\Gamma$.

Now let $B_{1},\dots,B_{n}$ be pairwise disjoint balls for local knots of $e$.  Since each $B_i$ meets $\Gamma$ in an arc in the interior of $e$, without loss of generality, we can assume that each $\partial B_i$ meets $N(\Gamma)$ in two disks contained in $N(e)-(N(v_1)\cup N(v_2))$.  For each $i$, let $F_i$ denote the annulus obtained from $\partial B_i$ by removing these two disks.   Since $N(e)\subseteq B$, for each $i$, $\partial F_i\subseteq B$.  By a standard cut-and-paste argument, we can assume that no components of $F_i-F$ are disks.  Thus if some $F_i$ is not contained in $B$, then there is an annulus component of $F_i-F$ in $\cl(M-B)$.  We can extend this annulus parallel to $F$ to obtain an annulus in $\cl(M-B)$ with boundaries in $\partial N(v_1)$ and $\partial N(v_2)$ which is incompressible in $\cl(M-B)$.  Hence, we can  obtain an isotopy of $(\cl(M-B),P')$ taking this annulus to $F$.  Thus we can remove any annuli from $F_i-F$ while pointwise fixing $\Gamma$ without introducing any new intersections of $F$ with some $F_j$.  This gives us an unknotting ball which contains $B_1\cup\dots\cup B_n$ as required.  \end{proof}

\medskip

\begin{lemma} Let $\Gamma$ be a 3-connected graph embedded in $S^3$ with locally knotted edges $e_1$,\dots, $e_n$.  Then there is a collection of pairwise disjoint unknotting balls for these edges.
\end{lemma}

\medskip\begin{proof}  By Theorem~1, we know that there is a collection $B_1$,\dots, $B_n$ of unknotting balls for $e_1$, \dots, $e_n$ respectively.  We see as follows that we can isotop $B_1$, \dots, $B_n$ to be pairwise disjoint while setwise fixing $\Gamma$.  First we will isotop $B_2$, \dots, $B_n$ off of $B_1$, then we will isotop $B_3$, \dots, $B_n$ off of $B_2$, and so on.

Let $C$ be a circle of intersection of $\partial B_1$ with $\partial B_2$ that bounds an innermost disk $\delta_1 $ on $\partial B_1$.  Since there are at least two innermost disks on $\partial B_1$ bounded by circles of $\partial B_1\cap\partial B_2$ and $e_1$ intersects $\partial B_1$ in precisely two points, we can choose the disk $\delta_1$ so that it intersects $e_1$ in at most one point.  Similarly, we can choose $\delta_2$ to be a disk bounded by $C$ in $\partial B_2$ which intersects $e_2$ in at most one point.  Since $B_i\cap e=\emptyset$ if $e\not =e_i$, the sphere $\delta_1\cup\delta_2$ meets $\Gamma$ in at most two points.  

By hypothesis, $\Gamma$ is 3-connected.  Thus a sphere cannot meet $\Gamma$ in precisely one point, and if a sphere meets $\Gamma$ in two points then those two points must be on the same edge.  Since the sphere $\delta_1\cup\delta_2$ meets each $e_i$ in at most one point, it follows that $\delta_1\cup\delta_2$ cannot intersect either $e_i$. Thus one component of $S^3-(\delta_1 \cup \delta_2)$ is a ball which is disjoint from $\Gamma$.  Using this ball, we can isotop $\delta_2$ to a disk parallel to $\delta_1$ by an isotopy of $B_2$ that fixes $\Gamma$.  This isotopy removes $C$ as a circle of intersection in $\bd B_1 \cap \bd B_2$.  Thus, by inducting on the number of circles of intersection, we can make $B_2$ disjoint from $B_1$.  Since we have changed $B_2$ by an isotopy pointwise fixing $\Gamma$, the new $B_2$ will still be an unknotting ball for $e_2$.   We continue this process to obtain pairwise disjoint unknotting balls for $e_1$,\dots, $e_n$.
 \end{proof}

\bigskip

\section{Adding local knots to spatial graphs}
In our study of topological symmetry groups, we would like to use local knots to prevent certain automorphisms from being induced by any homeomorphism of the embedding.  We begin with some more terminology.  

 By Theorem 1, we know that unknotting balls are unique up to an isotopy of $(S^3,\Gamma)$ fixing every vertex of $\Gamma$.   Thus if $B$ is an unknotting ball for an edge $e$, and the pair $(B, B\cap e)$ has knot type $K$, then we can unambiguously say that $e$ has {\it knot type $K$ in $\Gamma$} without making reference to a particular unknotting ball.  

Now let $B$ be a ball such that $B\cap \Gamma$ is an unknotted arc $A$ in an edge $e$.  Let $\Gamma'$ be obtained from $\Gamma$ by replacing $A$ with an arc $A'$ in $B$ such $(B, A')$ has knot type $K$ and $\partial A=\partial A'$.  Then we say that $\Gamma'$ is obtained from $\Gamma$ by {\it adding the local knot $K$ to $e$}.   

Suppose that $e$ has knot type $K_1$ in $\Gamma$.  Let $B$ be an unknotting ball for $e$, and let $B'\subseteq B$ such that $B'\cap e$ is an unknotted arc $A$.  Let $\Gamma'$ be obtained from $\Gamma$ by adding the local knot $K_2$ to $e$ within $B'$, and let $e'$ be the resulting edge.  Since $\Gamma-B=\Gamma'-B$, the ball $B$ is also an unknotting ball for $e'$ in $\Gamma'$.  Now by definition of the connected sum of knots, the pair $(B,B\cap e')$ has knot type $K_1\#K_2$.  Hence the edge $e'$ has knot type $K_1\#K_2$ in $\Gamma'$.

 \medskip

Finally, we will use the following terminology in the statement of the Knot Addition Lemma.  We say that an edge $e$ is {\it inverted} by an automorphism $\varphi$ if $\varphi$ interchanges the vertices of $e$.  If $X\subseteq Y$ is a set and $H$ is a group acting on $Y$, then we will use the notation $\langle X\rangle_H$ to denote the orbit of $X$ under $H$.

\medskip

\begin{add} \label{L:add}
Let $\Gamma$ be an embedding of a 3-connected graph in $S^3$ and let $H$ be a (possibly trivial) subgroup of $\TSG(\Gamma)$.  Let $e_1$, \dots, $e_n$ be edges of $\Gamma$ with distinct orbits under $H$.  Let $K_1$, \dots, $K_n$ be distinct prime knots, which are not local knots of $\Gamma$, such that $K_i$ is invertible if and only if $e_i$ is inverted by some element of $H$. Then for each $i$, the local knot $K_{i}$ can be added to the edges in $\langle e_i\rangle_H$ to create an embedding $\Gamma'$ such that $H \le \TSG(\Gamma')\leq \TSG(\Gamma)$.   Furthermore, for each $i$, let $e_i'$ be the embedding of $e_i$ in $\Gamma'$.  Then  $\langle e_i'\rangle_H=\langle e_i'\rangle_{\TSG(\Gamma')}$, and if $e_i'$ is inverted by an element of $\TSG(\Gamma')$ then $e_i'$ is also inverted by an element of $H$.
\end{add}

\medskip\begin{proof}   Let $X$ denote a set consisting of one point  in the interior of each edge of $\Gamma$.   For each $x\in X$, let $B_x$ denote a neighborhood of the point $x$ in the interior of $N(E)$ such that $B_x\cap \Gamma$ is an unknotted arc.  Let $S$ denote the set of all the $B_x$.  Thus the balls in $S$ are pairwise disjoint.  It follows from the uniqueness of $N(V)$, $N(E)$, and $X$ up to an isotopy fixing the vertices of $\Gamma$ that there is a group $G$ of orientation preserving diffeomorphisms of $(S^3,\Gamma)$ inducing $\TSG(\Gamma)$, such that for every $g\in G$, $g(N(V))=N(V)$, $g(X)=X$, and $g(S)=S$.  Note that $G$ need not be isomorphic to $\TSG(\Gamma)$.  Let $\widehat{H}$ denote a subgroup of $G$ which induces $H$ on $\Gamma$.  Since $S$ is setwise invariant under $G$, for any $B\in S$ we can define the orbit $\langle B\rangle_{\widehat{H}}\subseteq S$.

  For $i=1$, \dots, $n$, let $x_i$ denote the point of $X$ on $e_i$, and let $B_i=B_{x_i}$. Let $A_i$ be the arc $e_i\cap B_i$, let $A_i'$ denote an arc in $B_i$ with the same endpoints as $A_i$ containing the local knot $K_i$, and let $e_i'$ denote $e_i$ after $A_i$ has been replaced by $A_i'$.  By hypothesis, $K_i$ is prime and is not a local knot of $\Gamma$.  Also, $K_i$ is invertible if and only if $e_i$ is inverted by some element of $H$.  If $e_i$ is not inverted by any element of $H$, then we assign an orientation to $e_i$, which in turn induces an orientation on both $K_i$ and all of the edges of $\langle e_i\rangle_{H}$. Now let $e$ be an arbitrary edge of $\langle e_i\rangle_{H}$.  Then precisely one ball $B\in S$ intersects $e$.  Let $e'$ be obtained from $e$ by replacing the unknotted arc $B\cap e$ by an arc in $B$ with the same endpoints containing the local knot $K_i$ such that if $K_i$ is non-invertible then the orientation of $K_i$ with respect to the oriented edge $e'$ is the same as the orientation of $K_i$ with respect to the oriented edge $e_i'$.   If $e_i$ is inverted by some element of $H$, then $K_i$ was chosen to be invertible so the orientation of $K_i$ is not important.   Let $\Gamma'$ be obtained from $\Gamma$ by adding local knots to the edges in each orbit $\langle e_i\rangle_H$ in this way.   

 In order to prove that $H\leq \TSG(\Gamma')$, let $\varphi\in H$.  Since $H\leq \TSG(\Gamma)$, there is a diffeomorphism $g\in \widehat{H}$ which induces $\varphi $ on $\Gamma$.  We need to define an orientation preserving diffeomorphism of $(S^3,\Gamma')$ which induces $\varphi $ on $\Gamma'$.   For each $i$, let $\beta_i=\langle B_i\rangle_{\widehat{H}}$.  Then $\Gamma\cap (S^3-(\beta_1\cup\dots\cup \beta_n))=\Gamma'\cap (S^3-(\beta_1\cup\dots\cup \beta_n)$.  Hence we can define $g'|(S^3-(\beta_1\cup\dots\cup \beta_n))=g|(S^3-(\beta_1\cup\dots\cup \beta_n))$.  Also, since the set of vertices $V\subseteq (S^3-(\beta_1\cup\dots\cup \beta_n))$,  we know $g'|V=g|V$.

We extend $g'$ within each set of balls $\beta_i$ as follows.  Let $B\in \beta_i$, let $A$ denote the unknotted arc $B\cap \Gamma$, and let $A'$ denote the knotted arc $B\cap\Gamma'$.  First suppose that $g(B)\not = B$.  In this case, we let $A''$ denote the arc $g(B)\cap\Gamma'$.  Observe that since $g\in \widehat{H}$, the ball $g(B)\in\langle B_i\rangle_{\widehat{H}}$, and hence $A''$ contains the local knot $K_i$.   Thus the pairs $(B,A')$ and $(g(B),A'')$ both have knot type $K_i$, and if $K_i$ is non-invertible, then the knotted arcs $A'$ and $A''$ each has its orientation consistent with that of $e_i$.  Thus we can extend $g'$ within $B$ so that $g'((B, A'))=(g(B), A'')$.  Now suppose that $g(B)=B$.  If $g$ fixes the endpoints of $A$, then we can extend $g'$ to $(B,A')$ in such a way that $g'$ pointwise fixes the knotted arc $A'$.  If $g$ interchanges the endpoints of $A$ then $g$ inverts $e_i$ and hence the knot $K_i$ is invertible.  Thus we can extend $g'$ to $(B,A')$ in such a way that  $g'$ inverts $A'$. In this way we have extended $g'$ to every ball in $\beta_1\cup\dots\cup \beta_n$ such that $g'(\Gamma')=\Gamma'$.  Now  $g': (S^3,\Gamma')\rightarrow (S^3,\Gamma')$ induces $\varphi$ on $\Gamma'$.   It follows that $H\leq \TSG(\Gamma')$.

In order to prove that $\TSG(\Gamma')\leq \TSG(\Gamma)$, let $\varphi'\in\TSG(\Gamma')$.  Then $\varphi'$ is induced on $\Gamma'$ by an orientation preserving diffeomorphism $g'$ of $(S^3,\Gamma')$.  Since $H\leq \TSG(\Gamma')$, the orbit $E_i'=\langle e_i'\rangle_{H}$ is a set of edges in $\Gamma'$.  Suppose that each edge in $\langle e_i\rangle_{H}$ has knot type $J_i$ in $\Gamma$ (where $J_i$ might be the trivial knot).   Then each edge in $E_i'$ has knot type $J_i\#K_i$.  Since $\Gamma$ is 3-connected, it follows from \cite{FNT} that adding the local knot $K_i$ to an edge of $\Gamma$ does not cause any local knot to be added to any other edge of $\Gamma$.  Thus since $K_i$ is a prime knot that is not a local knot of $\Gamma$, the edges in $E_i'$ are the only edges in $\Gamma'$ containing $K_i$ among their local knots. It follows that for each $i$, $g'(E_i')=E_i'$.

By our construction, for each $e\in\langle e_i\rangle_{H}$, the neighborhood $N(e)$ is a ball for the local knot $K_i$ in the corresponding edge $e'$ of $\Gamma'$.  Thus by Theorem~1 and Lemma~1, we can choose a collection of pairwise disjoint unknotting balls for the edges in $E_1'\cup\dots\cup E_n'$ such that for each $e'\in E_i'$ the unknotting ball for $e'$ contains the ball $N(e)$.  For each $i$, let $\Delta_i$ denote the subset of these unknotting  balls which are unknotting balls for the edges in $E_i'$.  

Now $\Delta_i$ and $g'(\Delta_i)$ are each sets of unknotting balls for the edges of $\Gamma'$ in $E_i'$.  Since unknotting balls are unique up to isotopy by Theorem~1, there is an isotopy of $(S^3,\Gamma')$ fixing the vertices of $\Gamma'$ which takes $g'(\Delta_i)$ to $\Delta_i$.  Hence there is an orientation preserving diffeomorphism $f$ of $(S^3, \Gamma')$ fixing the vertices of $\Gamma'$ such that for each $i$, $f(g'(\Delta_i))=\Delta_i$.  Now $h=fg'$ is a diffeomorphism of $(S^3,\Gamma')$ which leaves each $\Delta_i$ setwise invariant and induces $\varphi'$ on $\Gamma'$.  For each $i$, the collection of balls in $\Delta_i$ contains both $\langle e_i\rangle_{H}$ and $E_i'=\langle e_i'\rangle_{H}$.  Thus $(S^3-(\Delta_1\cup\dots\cup \Delta_n))\cap \Gamma'=(S^3-(\Delta_1\cup\dots\cup \Delta_n))\cap \Gamma$.  So we can define $g|(S^3-(\Delta_1\cup\dots\cup \Delta_n))=h|(S^3-(\Delta_1\cup\dots\cup \Delta_n))$.  

We extend $g$ to the balls within each $\Delta_i$ as follows.  Let $B$ be one of the balls in $\Delta_i$, let $A_1$ denote the arc $B\cap \Gamma$ which is contained in some edge $e$ of $\langle e_i\rangle_{H}$, and let $A_1'$ denote the arc $B\cap \Gamma'$ which is contained in some edge $e'$ of $E_i'$.   Since $N(e)\subseteq B$ and $B$ is an unknotting ball for $e'$ in $\Gamma'$, $B$ must be an unknotting ball for $e$ in $\Gamma$ as well.  Thus since $e$ and $e'$ have knot types $J_i$ and $J_i\#K_i$ respectively, the pairs $(B, A_1)$ and $(B, A_1')$ also have knot types $J_i$ and  $J_i\#K_i$ respectively.

Suppose that $h(B)\not = B$.  Let $A_2=h(B)\cap\Gamma$ and $A_2'=h(B)\cap \Gamma'$.  Since $h((B,A_1'))=(h(B),A_2')$, the pair $(h(B),A_2')$ must also have knot type $J_i\#K_i$.  Since the edges in $E_i'$ are the only ones in $\Gamma'$ which contain $K_i$ among their local knots, $h(e')\in E_i'$.  Hence $h(e')$ is the embedding in $\Gamma'$ of some edge $\varepsilon \in\langle e_i\rangle_{H}$ in $\Gamma$.  Thus $\varepsilon$ has knot type $J_i$.  Since $h(B)\in \Delta_i$, it follows that the ball $N(\varepsilon)\subseteq h(B)$ and hence $h(B)$ is an unknotting ball for $\varepsilon$ in $\Gamma$.  Thus the pair $(h(B), A_2)$ has knot type $J_i$.  Now the pairs $(B,A_1)$ and $(h(B),A_2)$ both have knot type $J_i$.  Recall, that $K_i$ is a prime knot which is not contained in $\Gamma$.  In particular, $K_i$ is not among the prime factors of $J_i$.  Thus it follows from Schubert \cite{Sch} that if $J_i$ is non-invertible, then $J_i\#K_i$ is non-invertible as well.  Hence $h$ takes the oriented knot in $(B, A_1')$ to the oriented knot in $(h(B), A_2')$.  It follows that we can extend $g$ within $B$ so that $g((B,A_1))=(g(B), A_2)$.

Now suppose that $h(B)=B$.   If $h$ fixes the endpoints of $A_1'$, then we can extend $g$ to $(B,A_1)$ in such a way that $g$ pointwise fixes the arc $A_1$.  Suppose that $h$ interchanges the endpoints of $A_1'$.  Then $h$ inverts $(B,A_1')$.  Thus $J_i\#K_i$ must be invertible.  Since $J_i$ and $K_i$ have distinct knot types, it follows that $J_i$ is invertible.  Therefore, we can extend $g$ to $(B,A_1)$ in such a way that  $g$ inverts $A_1$. In this way we have extended $g$ to every ball in $\Delta_1\cup\dots\cup \Delta_n$ such that $g(\Gamma)=\Gamma$.  Since the vertices of $\Gamma$ are disjoint from $\Delta_1\cup\dots\cup \Delta_n$, the diffeomorphism $g: (S^3,\Gamma)\rightarrow (S^3,\Gamma)$ induces $\varphi'$ on $\Gamma$.   It follows that $\TSG(\Gamma')\leq \TSG(\Gamma)$.  

Finally, in order to show that $\langle e_i'\rangle_H=\langle e_i'\rangle_{ \TSG(\Gamma')}$, first observe that since $H\leq \TSG(\Gamma')$, $\langle e_i'\rangle_H\subseteq \langle e_i'\rangle_{ \TSG(\Gamma')}$.  Now let $e'\in \langle e_i'\rangle_{\TSG(\Gamma')}$.  Then for some $\varphi\in\TSG(\Gamma')$, we have $e'=\varphi(e_i')$.  Since $e_i'$ contains the local knot $K_i$, the edge $e'$ also contains the local knot $K_i$.  However, by our construction of $\Gamma'$, the only edges of $\Gamma'$ containing $K_i$ are the edges in $\langle e_i'\rangle_H$.  Thus $e'\in \langle e_i'\rangle_H$.  Hence $\langle e_i'\rangle_H=\langle e_i'\rangle_{ \TSG(\Gamma')}$.   Furthermore, suppose that $e_i'$ is inverted by some element of $\TSG(\Gamma')$.  Then $K_i$ must be invertible.  From our construction, it follows that $e_i$ must also be inverted by an element of $H$.  
\end{proof}

\medskip

The Finiteness Theorem below allows us to focus on topological symmetry groups which are induced by finite subgroups of $\Diff(S^3)$ (i.e., the group of orientation preserving diffeomorphisms of $S^3$). 

\medskip

\begin{finite} \cite{FNPT} \label{L:reembed}
Let $\Gamma$ be a 3-connected graph embedded in $S^3$.  Then there is an embedding $\Omega$ of $\Gamma$ in $S^3$ such that  $\TSG(\Gamma) \leq \TSG(\Omega)$ and $\TSG(\Omega)$ is induced by an isomorphic finite subgroup of $\Diff(S^3)$.
\end{finite}
\medskip

We use the Knot Addition Lemma together with the Finiteness Theorem to prove the following theorem.  Note that we use $\langle e_i\rangle_H$ to mean the orbit of the edge $e_i$ under the action of the group of automorphisms $H$.
\medskip

\begin{Subgroup thm}  Let $\Gamma$ be an embedding of a 3-connected graph in $S^3$ and let $H$ be a (possibly trivial) subgroup of $\TSG(\Gamma)$.  Let $e_1$, \dots $e_n$ be a set of edges in $\Gamma$ whose orbits under $H$ are distinct.  Suppose that any $\varphi\in \TSG(\Gamma)$ which  pointwise fixes $e_1$ and has $\varphi(\langle e_i\rangle_H)=\langle e_i\rangle_H$ for each $i$, also pointwise fixes a subgraph of $\Gamma$ that cannot be embedded in $S^1$.  Then there is an embedding $\Gamma'$ of $\Gamma$ with $H=\TSG(\Gamma')$.
\end{Subgroup thm}

\medskip\begin{proof}
 Let $K_1$, \dots, $K_n$ be distinct prime knots which are not local knots of $\Gamma$ and are invertible if and only if $e_i$ is inverted by some element of $H$.  We use the Knot Addition Lemma to add the local knot $K_i$ to every edge in $\langle e_i\rangle_H$ for each $i$.  Thus we obtain an embedding $\Gamma'$ such that $H\leq \TSG(\Gamma')\leq \TSG(\Gamma)$.  Furthermore, for each $i$, let $e_i'$ be the embedding of $e_i$ in $\Gamma'$.  Then  $\langle e_i'\rangle_H=\langle e_i'\rangle_{\TSG(\Gamma')}$, and if $e_i'$ is inverted by an element of $\TSG(\Gamma')$ then $e_i'$ is also inverted by an element of $H$.

In order to show that $\TSG(\Gamma')\leq H$, let $\alpha\in \TSG(\Gamma')$.  For each $i$, the edges in $\langle e_i'\rangle _H$ are the only edges of $\Gamma'$ containing the knot $K_i$.  Thus $\alpha(\langle e_i'\rangle _H)=\langle e_i'\rangle _H$.  Hence for some $h\in H$, $h(e_1')=\alpha(e_1')$.  Since $H\leq \TSG(\Gamma')$, $h^{-1}\alpha$ is an element of $\TSG(\Gamma')$ which setwise fixes $e_1'$.

If $h^{-1}\alpha$ inverts $e_1'$, then $e_1$ is inverted by some $f\in H$.  In this case, let $\varphi=fh^{-1}\alpha$.  Otherwise, let $\varphi=h^{-1}\alpha$.  In either case, $\varphi$ is an element of $\TSG(\Gamma')\leq \TSG(\Gamma)$ which pointwise fixes $e_1'$.  Also, for each $i$, $\varphi(\langle e_i'\rangle _H)=\langle e_i'\rangle _H$.  Thus $\varphi$ is an element of $\TSG(\Gamma)$ which pointwise fixes $e_1$ and  has $\varphi(\langle e_i\rangle_H)=\langle e_i\rangle_H$ for each $i$.  So by hypothesis, $\varphi$ pointwise fixes a subgraph of $\Gamma$ which cannot be embedded in $S^1$.

Now by the Finiteness Theorem there is an embedding $\Omega$ of $\Gamma$ in $S^3$ such that  $\TSG(\Gamma) \leq \TSG(\Omega)$ and $\TSG(\Omega)$ is induced by an isomorphic finite subgroup of $\Diff(S^3)$.  In particular, the automorphism $\varphi\in \TSG(\Omega)$ and hence is induced on $\Omega$ by a finite order $g\in\Diff(S^3)$.  Now $g$ pointwise fixes a subgraph of $\Omega$ which cannot be embedded in $S^1$.  Since $g$ has finite order, by Smith Theory \cite{Sm} $g$ must actually be the identity.  Thus $\varphi$ is the identity automorphism on $\Gamma'$.  Thus either $\alpha=hf^{-1}$ or $\alpha=h$.  In either case, $\alpha\in H$.  Hence $\TSG(\Gamma')=H$ as required.
\end{proof}
\medskip

Observe that the result below follows immediately from the Subgroup Theorem.

\medskip
\begin{Subgroup Corollary} \label{L:fundamentaledge}
Let $\Gamma$ be a 3-connected graph embedded in $S^3$ which has an edge $e$ that is not pointwise fixed by any non-trivial element of $\TSG(\Gamma)$.  Then for every (possibly trivial) subgroup $H$ of $\TSG(\Gamma)$, there is an embedding $\Gamma'$ of $\Gamma$ with $H = \TSG(\Gamma')$.
\end{Subgroup Corollary}

\bigskip

\section{Topological symmetry groups of complete graphs}

 In this section we use the Subgroup Corollary to prove Theorem~3.  We begin by stating some previous results that we will use.  The following theorem tells us which individual automorphisms of a complete graph $K_n$ can be induced by an orientation preserving diffeomorphism of some embedding of $K_n$ in $S^3$.  This result follows from Theorems 1 and 2 of \cite{Fl}.

\medskip

\begin{auto} \label{T:automorphism}
\cite{Fl} Let $n > 6$, and let $\varphi$ be an automorphism of $K_n$ of order $m$.  Then there is an embedding of $K_n$ in $S^3$ such that $\varphi $ is induced by an orientation-preserving diffeomorphism $h$ if and only if one of the following holds:
\begin{enumerate}
	\item $m$ is even and $m > 2$, all cycles of $\varphi $ are of order $m$, and $\varphi $ fixes no vertices;
	\item $m = 2$, all cycles of $\varphi $ are of order two, and $\varphi $ fixes at most two vertices;
	\item $m$ is odd, all cycles of $\varphi $ are of order $m$, and $\varphi$ fixes at most three vertices;
	\item $m$ is an odd multiple of 3, all cycles of $\varphi $ are of order $m$ except one of order 3, and $\varphi $ fixes no vertices.
\end{enumerate}
\end{auto}
\medskip

We will also use the following technical lemma, which follows from results in Section 2 of \cite{CFO}.  
\medskip

\begin{Z3}\label{Z3}\cite{CFO}  Let $\Omega$ be an embedding of $K_n$ in $S^3$ such that the group $\Z_3\times \Z_3$ is induced on $\Omega$ by a finite subgroup $H\leq\Diff(S^3)$.  Then there are at most two sets of 3 vertices which are each setwise invariant under $H$.  Furthermore, if a non-trivial $h\in H$ fixes any vertices of $\Omega$, then the set of fixed vertices of $h$ is equal to one of these sets of 3 vertices.  
\end{Z3}
\medskip

We will prove Theorem \ref{T:realize} by proving two propositions, one in which $\TSG(\Gamma)$ is cyclic or dihedral and the other in which $\TSG(\Gamma)$ is a subgroup of $D_m\times D_m$ for some odd $m$.  By the Subgroup Corollary, it suffices to prove in each proposition that there is an edge which is not pointwise fixed by any non-trivial element of $\TSG(\Gamma)$.

\medskip

\begin{prop} \label{P:rank1}
Let $\Gamma$ be an embedding of $K_n$  in $S^3$ with $n > 6$, such that $\TSG(\Gamma)$ is cyclic or dihedral.
 Then for any $H\leq \TSG(\Gamma)$, there is an embedding $\Gamma'$ of $K_n$  in $S^3$ such that $H = \TSG(\Gamma')$.
\end{prop}

\medskip\begin{proof} Suppose that $\TSG(\Gamma)$ is $\Z_m$ or $D_m$.  Then $\TSG(\Gamma)$ contains an element $\alpha$ of order $m$ which is induced by an orientation preserving diffeomorphism of $(S^3, \Gamma)$.  By the Automorphism Theorem, $\alpha$ has at least one $m$-cycle.  Let $v$ be a vertex in this $m$-cycle.  Then the edge $\overline{v\a(v)}$ is not pointwise fixed by any $\a^i$, with $1 \leq i \leq m-1$.  If $\TSG(\Gamma) = \Z_m$, then $\TSG(\Gamma)=\langle \alpha\rangle$.  Hence the result follows from the Subgroup Corollary.

Suppose that $\TSG(\Gamma) = D_m$.  Then there is an order 2 automorphism $\beta\in \TSG(\Gamma)$ which is induced by an orientation preserving diffeomorphism of $(S^3, \Gamma)$ such that $\a\b = \b\a^{-1}$.   Assume that for some $i$, $\beta\alpha^i$ pointwise fixes the edge $\overline{v\a(v)}$.  Then both $\beta\alpha^i(v)=v$ and $\beta\alpha^i(\a(v))=\a(v)$.  Thus $\a^i(v) = \b(v)$ since $\b$ has order 2.  Hence $\b \a^{i+1}(v) = \b \a \b(v) = \a^{-1}(v)$.  But we know that $\beta\alpha^i(\alpha(v))=\alpha(v)$.  Hence $\alpha^{-1}(v)=\alpha(v)$, and so $\a^2(v) = v$.  Since $v$ is in an $m$-cycle, this is impossible if $m > 2$.  Hence if $m>2$, then $\overline{v\a(v)}$ is not fixed by any element of $D_m$, and thus again the result follows from the Subgroup Corollary.

Now, we consider the case where $\TSG(\Gamma)=D_2$. In this case, $\a$ and $\b$ both have order 2 and together generate $D_2$.  By the Automorphism Theorem, $\a$, $\b$, and $\a\b$ each fix at most two vertices, and each have at least three 2-cycles since $n>6$.   So there are at least four vertices which are fixed by neither $\a$ nor $\b$.  If $\a$ and $\b$ agree on these four vertices, then the automorphism $\a^{-1}\b=\a\b$ would fix 4 vertices.  As this is contrary to the Automorphism Theorem, there is at least one vertex $v$ which is fixed by neither $\a$ nor $\b$ such that $\a(v) \neq \b(v)$.  

It follows that $v$, $\a(v$), $\b(v)$ are all distinct vertices.  We see that $\a\b(v)$ is also distinct from these three vertices as follows.  If $\a\b(v) = v$, then $\b(v) = \a(v)$; if $\a\b(v) = \a(v)$, then $\b(v) = v$; and if $\a\b(v) = \b(v)$ then $\a(v) = v$.  All three cases are impossible since $v$, $\a(v$), $\b(v)$ are distinct.  Hence $v$, $\a(v$), $\b(v)$, and $\a\b(v)$ are all distinct.  It follows that none of $\a$, $\b$ or $\a\b$ pointwise fixes $\overline{v\a(v)}$.  Now our result follows again from the Subgroup Corollary.
\end{proof}

\medskip

The remaining cases in the proof of Theorem \ref{T:realize} are when $\TSG(\Gamma)$ is a subgroup of $D_m \times D_m$ for some odd $m$.   These subgroups are described by the following lemma.  This result is undoubtedly known, but we were unable to find a reference, so we include a proof here.
\medskip

\begin{lemma} \label{T:subgroups}  Let $m\geq 3$ be odd, and let $G$ be a non-trivial subgroup of $D_m \x D_m$.  Then $G$ is isomorphic to
one of the following groups where $r$, $s\geq 3$ are odd :
$\Z_2$,
$\Z_r$, $\Z_{2r}$, $D_2$,
$D_r$, $D_{2r}$,
$\Z_r \x \Z_s$,
$D_r \x \Z_s$,
$D_r \x D_s$,
or
$(\Z_r \x \Z_s) \rtimes \Z_2$ such that for any nontrivial elements $g\in\Z_r \x \Z_s$ and $\varphi \in\Z_2$ we have $g\varphi=\varphi g^{-1}$.\end{lemma}

\medskip\begin{proof}
Each element of $D_m$ can be thought of as either a rotation or a reflection of a circle.
Define a homomorphism
$\rho: D_m \to \Z_2$ by
$\rho(x) = 1$ if $x$ is a reflection and
$\rho(x)=0$ otherwise.
Let $\sigma = (\rho \x \rho)|_G : G \to \Z_2 \x \Z_2$.

Now, $H = \kernel(\sigma) \sbgp \kernel(\rho \x \rho) = \Z_m \x \Z_m$.
Every subgroup of $\Z_m \x \Z_m$ has rank at most 2, and hence $\rank(H)\leq 2$.
Thus $H$ is either trivial or isomorphic to
$\Z_r$ or $\Z_r \x \Z_s$ for some odd $r, s \ge 3$.  If $H$ is trivial,
then $\sigma$ is an isomorphism
from $G$ to a subgroup of $\Z_2 \x \Z_2$.
Thus $G$ is either the trivial group, $\Z_2$, or $D_2=\Z_2 \x \Z_2$, so we are done.  Hence we shall assume that $\rank(H) = 1$ or 2.

Let $k =\rank(\sigma(G))$, then $k \leq 2$.  If $k = 0$, then $G = H$, and we are done.  Thus we shall assume that  $k =1$ or 2.   If $k=1$ then without loss of generality we can assume that the generator of $\sigma(G)$ is either $(1,1)$ or $(1,0)$, and if $k=2$ then we can assume the generators of $\sigma(G)$ are $(1,0)$ and $(0,1)$.  We  consider two cases, according to whether or not $\sigma(G)$ is generated by $(1,1)$.
\medskip

\noindent {\bf Case 1:} $\sigma(G)$ is generated by $(1,1)$.

 In this case, $G=\langle H,(x,y)\rangle$ where both $x$ and $y$ are reflections.  Now for every nontrivial $h\in H$,  $hx=xh^{-1}$ and $hy=yh^{-1}$.  Thus, $G$ is either isomorphic to $D_r$ or $(\Z_r \x \Z_s) \rtimes \Z_2$ depending on whether $H$ is isomorphic to
$\Z_r$ or $\Z_r \x \Z_s$, respectively.  Thus, in Case 1 we are done.  

\medskip

Before we consider Case 2 we make the following observation.  Note here we use $e$ to denote the identity element.
\medskip

\noindent\textbf{Observation}:  If $(x,y) \in G$ where
$x$ is a reflection and  $y$ a rotation, then $(e,y) \in H$ and $(x,e)\in G$.  If $(x,y) \in G$ where
$x$ is a rotation and  $y$ a reflection, then $(x,e) \in H$ and $(e,y)\in G$.

\medskip

\noindent\textbf{Proof of Observation}: We prove the first assertion as follows.  We know that $x^{m+1} = e$, since $m$ is odd and $x$ is a reflection. Also, $y^m =e$ since $y \in \Z_m \le D_m$.
Thus $(x,y)^{m+1}=(e,y) \in G$.  Now since $y$ is a rotation,
it follows that $(e,y) \in H$.  Also, since $(x,y)$ and $(e,y)$ are both in $G$, we have $(x,e)\in G$.  The second assertion is proved similarly.

\medskip
We will use this Observation in the proof of Case 2.

\medskip

\noindent {\bf Case 2:} $\sigma(G)$ is not generated by $(1,1)$.

Without loss of generality we can assume that either $\sigma(G)=\langle(1,0)\rangle$  or $\sigma(G)=\langle(1,0),(0,1)\rangle$.  In the first case, $G=\langle H,(x_1,y_1)\rangle$, where
$x_1$ is a reflection and  $y_1$ is a rotation; and in the second case, $G=\langle H,(x_1,y_1), (x_2,y_2)\rangle$, where $x_1$ and $y_2$ are  reflections and  $y_1$ and
$x_2$ are rotations.  Now, by applying the Observation we see that in the first case $(e,y_1) \in H$ and $(x_1,e)\in G$, and   
in second case, $(e,y_1),(x_2,e) \in H$ and $(x_1,e), (e,y_2)\in G$.  Thus either $G=\langle H,(x_1,e)\rangle$ or $G=\langle H,(x_1,e)(e,y_2)\rangle$.

In either case, $(x_1,e)\in G$ and $x_1$ is a reflection.  Thus for any $(a,b)\in H$, we know that $(ax_1,b)\in G$, and $ax_1$ is a reflection and $b$ is a rotation.  So we can apply the Observation to conclude that $(e,b)\in H$.  Now since $(e,b)$ and $(a,b)$ are both in $H$, it follows that $(a,e)\in H$.
\medskip

\noindent {\bf Subcase (a)} $\rank(H)=1$.

 In this case, without loss of generality we can assume that $H=\langle(a,b)\rangle$, where $a$ and $b$ are rotations and $(a,b)$ is of odd order $r$.  Then either $G=\langle(a,b),(x_1,e)\rangle$ or $G=\langle(a,b),(x_1,e), (e,y_2)\rangle$.   As we saw above $(e,b)\in H$.  Hence  for some $j$, $(e,b) = (a,b)^j$.
This gives $a^j = e$ and $b^{j-1} = e$.  Assume that neither $a$ nor $b$ is the identity.  Then $j$ and $j-1$ must each divide $m$ or equal 0, since $H\leq \mathbb{Z}_m\times \mathbb{Z}_m$.
Since $m$ is odd, $j$ and $j-1$ cannot both divide $m$, unless $j=1$.  Thus either $j=0$ or $j=1$.  If $j=0$, then $b=e$ and if $j=1$ then $a=e$.  Hence without loss of generality we assume that $a=e$.

Now $G=\langle(e,b),(x_1,e)\rangle$ or $G=\langle(e,b),(x_1,e), (e,y_2)\rangle$.  In the former case, $G\cong \Z_2\times \Z_r\cong \Z_{2r}$, since $r$ is odd.  In the latter case, $G\cong \Z_2\times  D_r$, since $(e,b)$ and $(e,y_2)$ generate $D_r$
and they both commute with $(x_1,e)$.  Now since $r$ is odd, $G\cong \Z_2\times  D_r\cong D_{2r}$.

\medskip

\noindent {\bf Subcase (b)} $\rank(H)=2$.

In this case, without loss of generality we can assume that $H=\langle(a,b),(c,d)\rangle$, where $a$, $b$, $c$, $d$, are rotations such that $(a,b)$ and $(c,d)$ are of odd orders $r$ and $s$ respectively.  Hence either $G=\langle(a,b),(c,d),(x_1,e)\rangle$ or $G=\langle(a,b),(c,d),(x_1,e), (e,y_2)\rangle$.  In either case, by arguing as we did before Subcase (a), we conclude that $(a,e)$, $(e,b)$, $(c,e)$, $(e,d)\in H$.  Thus $H = \langle (a,e),(e,b),(c,e),(e,d) \rangle$.  It follows that for some $u$ and $v$ of order $r'$ and $s'$ respectively, $H = \langle (u,e),(e,v) \rangle$.  Also $r'$, $s' \ge 3$ and are odd.  Thus either $G = \langle (u,e),(e,v), (x_1,e) \rangle$ or $G = \langle (u,e),(e,v), (x_1,e), (e,y_2) \rangle$.  In the former case, $G\cong D_{r'} \x \Z_{s'}$, and in the latter case $G \cong D_{r'} \x D_{s'}$.  Thus again we are done.\end{proof}

\medskip

\begin{prop} \label{P:product}
Let $\Gamma$ be an embedding of $K_n$ in $S^3$ with $n > 6$, such that $\TSG(\Gamma)$ is a subgroup of $D_m\times D_m$ for some odd $m$.  Then for every (possibly trivial) subgroup $H$ of $\TSG(\Gamma)$, there is an embedding $\Gamma'$ of $K_n$  in $S^3$ such that $H = \TSG(\Gamma')$.
\end{prop}

\medskip\begin{proof}  In Proposition~1 we have dealt with the case where $\TSG(\Gamma)$ is cyclic or dihedral.  So we shall assume that $\TSG(\Gamma)$ is neither.  Now it follows from Lemma \ref{T:subgroups} that $\TSG(\Gamma)$ is isomorphic to one of the following groups, where $r$, $s\geq 3$ are odd:
\begin{enumerate}
\item $\Z_r \times \Z_s$
	\item $D_r \times \Z_s$ 
	\item $D_r \times D_s$ 
	\item $(\Z_r \times \Z_s) \ltimes \Z_2$ such that for any nontrivial elements $g\in\Z_r \x \Z_s$ and $\varphi \in\Z_2$ we have $g\varphi=\varphi g^{-1}$\end{enumerate}

Observe that in all 4 cases, $\Z_r \times \Z_s  \leq \TSG(\Gamma)$.  We choose $\alpha$, $\beta\in  \TSG(\Gamma)$ such that $\langle\alpha,\beta\rangle=\mathbb{Z}_r\times \mathbb{Z}_s$, the order of $\a$ is $r$, the order of $\b$ is $s$, and $\langle\alpha\rangle\cap \langle\beta\rangle$ only contains the identity.  Recall that for any vertex $v$, $\langle v\rangle_{\alpha}$ and $\langle v\rangle_{\beta}$ denote the $\alpha$-orbit and $\beta$-orbit respectively of $v$.  
 
 We prove in the following two cases that there is a vertex $v$  which is contained in both an $r$-cycle of $\alpha$ and an $s$-cycle of $\beta$ such that the two cycles intersect only in $v$.
 \medskip
 
\noindent {\bf Case 1:}  Not both $r=3$ and $s=3$.  

By the  Automorphism Theorem, there are at most 3 vertices which are not contained in an $r$-cycle of $\alpha$ and at most 3 vertices which are not contained in an $s$-cycle of $\beta$.  Since $n>6$, we can pick a vertex $v$ which is contained in both an $r$-cycle of $\alpha$ and an $s$-cycle of $\beta$.   Now suppose that for some $i$ and $j$ such that $0<i<r$ and $0<j<s$, $\beta^j(v)=\alpha^i(v)$.  Let $w\in \langle v\rangle_{\alpha}$, then $w = \alpha^k(v)$ for some $k < r$.  Hence $\alpha^i(w) = \alpha^{i+k}(v) = \alpha^k(\beta^j(v)) = \beta^j(\alpha^k(v)) =\beta^j(w)$.  Thus $\alpha^i\beta^{-j}$ fixes $w$, and hence fixes every vertex in $\langle v\rangle_{\alpha}$.  Similarly, $\alpha^i\beta^{-j}$ fixes every vertex in $\langle v\rangle_{\beta}$. We know that $\alpha^i\beta^{-j}$ is not the identity since  $\langle \alpha\rangle\cap \langle \beta\rangle$ only contains the identity.  Thus by the Automorphism Theorem, $\alpha^i\beta^{-j}$ fixes at most 3 vertices.   Therefore $\langle v\rangle_{\alpha}\cup \langle v\rangle_{\beta}$ contains at most 3 vertices.  On the other hand, we chose $v$ so that $|\langle v\rangle_{\alpha}| =r$ and $|\langle v\rangle_{\beta}|=s$. Thus $r\leq 3$ and $s\leq 3$.   This is a contradiction since  $r$, $s\geq 3$ and they do not both equal 3. Thus we conclude that there is a vertex $v$ which is contained in both an $r$-cycle of $\alpha$ and an $s$-cycle of $\beta$ such that $\langle v \rangle_{\a} \cap \langle v \rangle_{\b} = \{v\}$.

\medskip
 
\noindent {\bf Case 2:} $r=s=3$.

By the Finiteness Theorem, there is a re-embedding $\Omega$ of $\Gamma$ such that $\TSG(\Gamma) \leq \TSG(\Omega)$ and $\TSG(\Omega)$ is induced by an isomorphic subgroup of $\Diff(S^3)$.  Thus $\langle \alpha,\beta\rangle=\Z_3\times \Z_3$ is induced on $\Omega$ by a finite subgroup of $\Diff(S^3)$.  Now it follows from the $\mathbb{Z}_3\times\mathbb{Z}_3$ Lemma that there are at most two sets of 3 vertices which are each setwise invariant under both $\alpha$ and $\beta$, and if either $\alpha$ or $\beta$ has any fixed vertices then the set of its fixed vertices is equal to one of these sets of 3 vertices.  Since $n>6$, we can pick $v$ to be a vertex which is not in one of the above sets of 3 vertices. Thus $\langle v\rangle_{\alpha}$ contains 3 vertices and cannot be setwise fixed by $\beta$. Now $|\langle v\rangle_{\beta}|=3$  and $\langle v\rangle_{\alpha}\not=\langle v\rangle_{\beta}$.  Suppose that for some $i$ and $j$, $0<i<3$ and $0<j<3$, we have $\beta^j(v)=\alpha^i(v)$.  Since $r=s=3$, $\alpha^i=\alpha^{\pm1}$ and $\beta^j=\beta^{\pm1}$.  Thus $\alpha(v)=\beta^{\pm1}(v)$.  Thus $\langle v\rangle_{\alpha}=\langle v\rangle_{\beta}$, which is contrary to our choice of $v$.  Hence there is a vertex $v$ which is contained in both an $r$-cycle of $\alpha$ and an $s$-cycle of $\beta$ such that $\langle v \rangle_{\a} \cap \langle v \rangle_{\b} = \{v\}$.

\medskip

Let $v$ be the vertex given by Case 1 or 2, and let $e$ be the edge $\overline{v\a(\b(v))}$.  Now $\a^i \b^j(v)\not =v$ for all $0\le i<r$ and $0\le j<s$ except when both $i= 0$ and $j = 0$.  Thus the edge $e$ is not pointwise fixed by any non-trivial element of $\Z_r \times \Z_s$.  If $\TSG(\Gamma)=\Z_r \times \Z_s$, then we can apply the Subgroup Corollary to conclude that for every subgroup $H\leq \TSG(\Gamma)$, there is an embedding $\Gamma'$ of $K_n$ in $S^3$ such that $H = \TSG(\Gamma')$.

Suppose that $\TSG(\Gamma) = D_r \times \Z_s$.  Then there is an order 2 automorphism $\g\in \TSG(\Gamma)$ such that $\alpha\g = \g\alpha^{-1}$ and $\beta\g= \g\beta$.  We saw above that the edge $e$ is not pointwise fixed by any non-trivial element of the subgroup $\langle \alpha,\beta\rangle$ of $\TSG(\Gamma)$.  Assume for the sake of contradiction that some  $\g \a^i \b^j$ pointwise fixes $e$.  This means that both $\g \a^i \b^j (v) = v$ and $\g \a^i \b^j (\a\b(v)) = \a\b(v)$.  So $\a^i\b^j(v) = \g(v)$, and hence:
$$\a^{-1}\b(v)=\g^2 \a^{-1}\b(v)=\g\a\b(\g(v))=\g\a\b(\a^i\b^j(v))=\g\a^i\b^j(\a\b(v))=\a\b(v)$$  
 
 \noindent Thus $\a^2\beta(v) =\beta(v)$, and hence $\a^2(v)=v$.  But $v$ is contained in an $r$-cycle of $\a$ and $r$ is odd.  So this is impossible.  Therefore $e$ is not pointwise fixed by any non-trivial element of $D_r \times \Z_s$.  Hence again we can apply the Subgroup Corollary to conclude that for every subgroup $H\leq \TSG(\Gamma)$, there is an embedding $\Gamma'$ of $K_n$ in $S^3$ such that $H = \TSG(\Gamma')$.

Next suppose that $\TSG(\Gamma) = D_r \times D_s$.  Then, in addition to $\a$, $\b$, and $\g$, there is an order 2 automorphism $\d\in \TSG(\Gamma)$ such that $\alpha\d = \d\alpha$, $\beta\d= \d\beta^{-1}$, and $\d\g=\g\d$.  We saw above that the edge $e$ is not pointwise fixed by any non-trivial element of the form $\a^i\b^j$ or $\g\a^i\b^j$.  By an analogous argument we see that $e$ is also not pointwise fixed by any non-trivial element of the form $\d\a^i\b^j$.  Now for the sake of contradiction assume that for some $i$ and $j$, $\g\d\a^i\b^j$ pointwise fixes $e$.  Thus both $\g\d\a^i\b^j(v) = v$ and $\g\d\a^i\b^j(\a\b(v)) = \a\b(v)$.  So $\a^i\b^j(v) = \g\d(v)$, and hence:
 $$\a^{-1}\b^{-1}(v)=(\g\d)^2 \a^{-1}\b^{-1}(v)=\g\d\a\b(\g\d(v))=\g\d\a\b(\a^i\b^j(v))=$$
 $$\g\d\a^i\b^j(\a\b(v))=\a\b(v)$$

 \noindent  Thus $\a^{r-2}(v) =\a^{-2}(v)=\beta^2(v)$.  However, since $r$, $s\geq 3$, we have contradicted the fact that $\langle v \rangle_{\a} \cap \langle v \rangle_{\b} = \{v\}$.   So, once again, the edge $e$ is not fixed by any non-trivial element of $\TSG(\Gamma)$, and thus we can again apply the Subgroup Corollary.

Finally, suppose that $\TSG(\Gamma) = (\Z_r \times \Z_s) \ltimes \Z_2$.  Then in addition to $\a$ and $\b$, $\TSG(\Gamma)$ contains an order 2 automorphism $\mu$ such that $\a\mu = \mu\a^{-1}$ and $\b\mu = \mu\b^{-1}$.   If for some $i$ and $j$, $\mu\a^i\b^j$ pointwise fixes $e$, then by substituting $\mu$ for $\g\d$ in the previous paragraph we would again obtain a contradiction.  Thus yet again we conclude that $e$ is not fixed by any non-trivial element of $\TSG(\Gamma)$, so we can apply the Subgroup Corollary.
\end{proof}

\medskip

Theorem \ref{T:realize} now follows immediately from Propositions~1 and 2.
\bigskip

\section{Acknowledgments}  The authors would like to thank Curtis Bennett and Michael Aschbacher for several very helpful conversations which led to significant improvements in this paper.  We also would like to thank the anonymous referee for helpful suggestions.

\end{document}